\documentclass[a4paper,12pt]{article}
\usepackage[utf8]{inputenc}
\usepackage[english]{babel}
\usepackage{amsmath,amsfonts,amssymb} 
\usepackage{graphicx}
\usepackage{enumerate}
\usepackage{url}

\pdfoutput=1

\usepackage{special2_iiets_thue-morse}

\newcommand{\theofont}{\sf \normalsize}
\newtheorem{Def}{\theofont Definition}

\newtheorem{Lem}{\theofont Lemma}
\newtheorem{Prop}{\theofont Proposition}

\newtheorem{Theo}{\theofont Theorem}
\newtheorem{Cor}{\theofont Corollary}
\newcommand{\PR}{{\vspace*{-0.1cm}\theofont{\em{Proof.}} } }

\newcommand{\EPR}{\hfill $\Box$ \vspace*{0.4cm}}
\newcommand{\EEX}{\hfill $\Diamond$ \end{Ex} \vspace*{0.3cm}}
\newcommand{\EQ}{\begin{equation}}
\newcommand{\EEQ}{\end{equation}}
\newcommand{\smbf} {\footnotesize \bf}
\newcommand{\II} {[0,1)}

\newcommand{\DiscPre}{{\cal D}}%
\newcommand{\Disc}{\overline{\cal D}}
\newcommand{\Discaccr}{{\cal D}_{acc,r}}%
\newcommand{\Discaccl}{{\cal D}_{acc,l}}%
\newcommand{\Lprim}{L'}
\newcommand{\Lmprim}{L_m'}
\newcommand{\LprimOne}{L_\diamond'}
\newcommand{\phiinv}{\phi_\mu^{inv}}
\newcommand{\phiLprim}{\widehat{\phi}_\mu}
\newcommand{\Ninj}{Ninj_{T_L}}
\newcommand{\Ninjm}{Ninj_{T_{L_m}}}
\newcommand{\ConsecL}{Consec_L}
\newcommand{\ConsecLm}{Consec_{L_m}}
\newcommand{\ConsecPairL}{ConsecPair_L}
\newcommand{\ConsecPairLm}{ConsecPair_{L_m}}
\newcommand{\ShiftsAM}{Shifts_{am}}
\newcommand{\ShiftsAMZ}{Shifts_{amz}}

\begin{document}

\title{
  {\large \bf \sffamily  Non-Injectivity of \\ \vspace*{0.1cm}
    Infinite Interval Exchange Transformations\\ \vspace*{-0.12cm}
    and Generalized Thue-Morse Sequences}
}
\date{}

\author{
  \small Luis-Miguel Lopez \\
{\scriptsize Tokyo University of Social Welfare, Japan}
  \and
  \small Philippe Narbel\\
{\scriptsize LaBRI, University of Bordeaux, France}
}

\maketitle

\vspace*{-0.44cm}
\noindent
    {\footnotesize {\sc Abstract.} {\em In this paper we study the non-injectivity
        arising in infinite interval exchange transformations. In particular, we
        build and analyze an infinite family of infinite interval exchanges
        semi-conjugated to generalized Thue-Morse subshifts, whose
        non-injectivity occurs at a characterizable finite set of points.}}

\vspace*{0.3cm}

\noindent
{\scriptsize {\sc MSC2020 Mathematics Subject Classification:} 37E05, 37B10,
  37A05.}

%


%
\section{Introduction}
Interval exchange transformations are functions over $\II$
which induce archetypal measure-preserving one-dimensional dynamical systems
coming from first-return maps of measure-preserving flows on surfaces
\cite{Ose66,Kea75,Man87,KH95}.
In the classic case, such maps are defined as bijective piecewise isometries,
each determined by a finite set of discontinuities
and a permutation of their continuity intervals.
When infinitely many discontinuities are allowed,
the basic properties of these maps --~local isometry and measure
preservation~-- are straightforwardly kept,
but not bijectivity.
Of course when non-bijectivity occurs on a null measure set, many general
dynamical properties are not affected. Nevertheless when it comes to
non-injectivity, it implies non-invertibility and backward bifurcations in the
corresponding flows, thus a significant structural effect on the overall
dynamics.

In this paper we study the emergence of this non-bijectivity
by analyzing infinite interval exchange transformations showing a
characterizable finite set of non-injectivity points.
These interval exchanges are obtained as semi-conjugates of classic symbolic
dynamical systems, i.e., {\em subshifts},
 induced by a generalization of the classic {\em Thue-Morse
   sequence}~\cite{Pro1851,Thu12,Mor21b} over $m \geq 2$
 letters~\cite{TS95,AS03}
using a technique developed in~\cite{LN17}, 
so that the non-injectivity of the former can be studied through the
non-injectivity of the latter using symbolic dynamic
techniques.
As a main result, we exhibit an infinite family of infinite interval exchange
transformations indexed by $m \geq 2$
whose non-injectivity occurs exactly on $m(m-1)$ points sent onto $m-1$ points
(see Theorem~\ref{non_inj_th}, p.~\pageref{non_inj_th}). 
\section{Infinite Interval Exchange Transformations}
\label{iiets_sec}

A {\em finite interval exchange transformation} $T$ is a measure-preserving map
which occurs as a bijective piecewise isometry of the semi-open unit interval
$\II$, having finitely many discontinuities, generally taken~as right-continuous
and orientation-preserving.  Its finite set $\DiscPre$ of discontinuities
determines a finite set of continuity intervals for $T$ on which local
isometries act globally as a permutation (an ``exchange'') of their order
in $\II$.

When an interval exchange $T$ is {\em infinite}, i.e., it has an infinite set
$\DiscPre$ of discontinuities,
every property of the finite case can be retained
except bijectivity because of what may happen at accumulation points of
discontinuities.
Since infinity does not necessarily preclude bijectivity, in particular
injectivity, many of the existing definitions of infinite exchanges just require
it to hold.
However, in a full setting, 
a definition of such a map~$T$ allowing non-bijectivity is the following:

\begin{Def}\label{iiet_def}
An orientation-preserving {\smbf interval exchange transformation} is a right
continuous and measure-pre\-serving map $T: \II \rightarrow \II$
with a countable set $\DiscPre$ of discontinuities such that: 

\vspace*{-0.2cm}
 \begin{itemize}\itemsep=-0.05cm
 \item The closure~$\Disc$ of $\DiscPre$ in $\II$
 has null measure, and determines a partition of~$\II$ formed by:

 \vspace*{-0.2cm}
 \begin{enumerate}[i.]\itemsep=0cm
   \item The right-open intervals $[x, x')$ with $x \in \Disc$, $x' \in \Disc
    \cup \{1\}$, so that $[x, x') \cap \Disc = \{x\}$ (these intervals are
      maximal with regard to $\Disc$);
    
   \item The set $\Discaccr$ of accumulation points of $\Disc$ from the right in
     $\II$.
 \end{enumerate}

\item $T$ is a translation on each right-open interval of the above partition.
\end{itemize}
\end{Def}

\noindent
With regard to bijectivity, using that $T$ is measure-preserving and defined as
a set of translation maps we first classically have:
\begin{Prop}\label{non-inj_first_prop}
  $T$ is injective on $\II \setminus \Disc$.
\end{Prop}
\noindent
Then in the finite case,
Proposition~\ref{non-inj_first_prop} extends well~\cite{Man87}:
\begin{Prop}\label{disc_finite_prop}
  If $\Disc$ is finite, then $T$ is bijective on $\II$.
\end{Prop}
\noindent
In the infinite case, injectivity and surjectivity do not necessarily hold.
A first situation is the following: 
\begin{Prop}\label{discaccr_empty_prop}
  If $\Discaccr = \emptyset$, then $T$ is injective on $\II$.
\end{Prop}
\PR Assume there exist $x,x' \in \Disc$ such that $T(x)=T(x')$. Since
$\Discaccr$ is empty, $x$ and $x'$ must be the left ends of distinct
continuity intervals. But then non-injectivity would hold on 
subintervals starting at $x$ and $x'$, impairing measure preservation.
\EPR

\vspace*{0.1cm}
\noindent
When $\Disc$ is infinite, the above case
implies that the set $\Discaccl$ of accumulation points of $\Disc$ from the left
in $\II$ is not empty. But accordingly, their image limit points do not
determine images for~$T$. This fact explains that:
\begin{Prop}\label{non_surj_prop}
  If $\Disc$ is infinite with $\Discaccr = \emptyset$, then $T$ can be
  non-surjective.
\end{Prop}

\noindent
In fact, many occurrences of infinite interval exchanges in the literature
happen to be of the above kind.
For instance, the {\em von Neumann-Kakutani transformations} over
$\II$ (also called {\em van der Corput maps})~\cite{Fri70,Fog02}
are defined for each $b \geq 2$ as $f(x)=x-1+b^{-n}+b^{-n-1}$ for $x
\in [1-b^{-n},1-b^{-n-1})$, with $n \in \mathbb{N}$.
They have a unique accumulation point at~$1$ from the left, they are injective
and their ranges do not include $0$ (see the graph for $b=2$ at the left of the
figure below).
The same~holds for the maps obtained
in~\cite{AOW85} and the rational one in~\cite{HRR20}.
The derived {\em Kakutani-Fibonacci transformation}~\cite{CJV14} is also
injective and non-surjective, with two accumulation points from the left, one of
them in $(0,1)$
(see its graph below at the right):%

\vspace*{-0.1cm}
\begin{center} 
  \includegraphics[width=4cm]{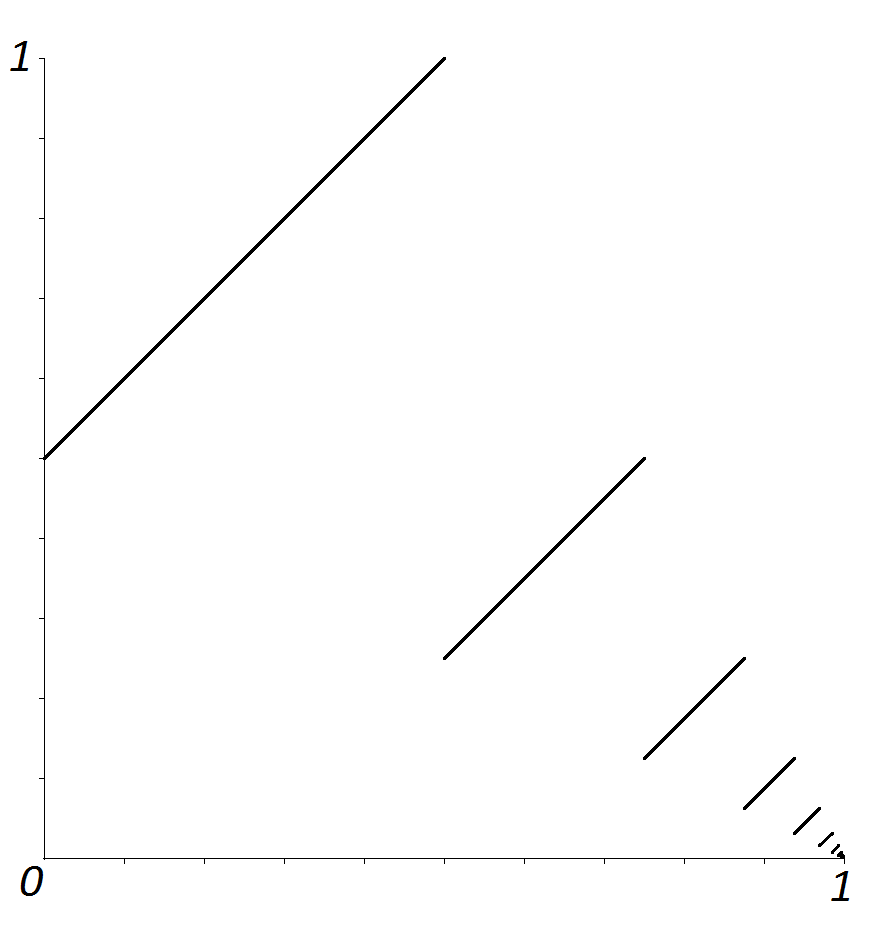}
  \hspace*{1.5cm}
  \includegraphics[width=4cm]{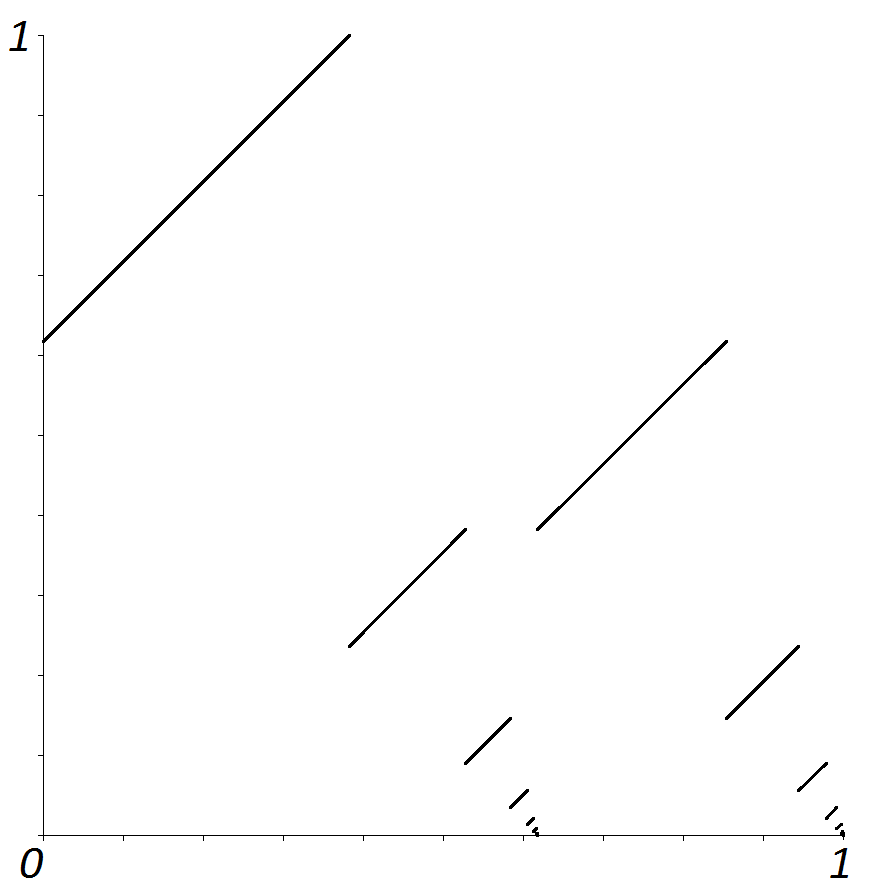}
\end{center}

Now, when $\Discaccr \neq \emptyset$, non-injectivity may occur.
A first property is that~$\Discaccr$ takes a primary role in non-injectivity
with regard to $\Disc$:
\begin{Prop}\label{discaccr_overwhelm_prop}
If $T$ is non-injective on $x \in \Disc \setminus \Discaccr$, then $U_x=\{ x'
\in \II \mid T(x') = T(x), \; x'\neq x \}$ has all its points in $\Discaccr$.
\end{Prop}
\PR If $x' \in U_x$, then for the same reason as in the proof of
Proposition~\ref{discaccr_empty_prop}, $x'$ cannot belong to an interval of
continuity of $T$.
\EPR

\noindent
Then, a case where non-injectivity indeed occurs when $\Discaccr \neq
\emptyset$ is for infinite interval exchanges which are {\em piecewise
  increasing on a finite partition}
(recall that {\em piecewise} means to hold on a locally finite family of
specific subsets --~e.g., the above von Neumann-Kakutani maps are piecewise
increasing but not on a finite partition, whereas the examples in
Section~\ref{ex_sec} do):

\begin{Prop}\label{piecewise_non-inj_prop}
  Let $T$ be piecewise increasing on a finite partition and such that $\Discaccr
  \neq \emptyset$.  Then $T$ is non-injective on $\Discaccr$, and only on this
  set.
\end{Prop}
\PR Let $x \in \Discaccr$. Then there is an interval $[x, y)$ where $T$ is
  increasing, and where $[x, y)$ is a union of infinitely many intervals with
    left endpoints in $\DiscPre$.
Let us consider the set of their images by $T$. Since~$\DiscPre$ is made
of discontinuities of $T$, there are infinitely many intervals in between these
images whose left endpoints accumulate at $T(x)$.
The inverse images by $T$
form at least one accumulation point $x'$, also from the right since $T$ is
piecewise increasing on a finite partition, so that $x' \in \Discaccr$ too.
By right continuity, $T(x') = T(x)$ with $x' \neq x$ since $x'$ is an
accumulation point of left interval endpoints lying outside $[x, y)$.  Hence,
  $T$ is non-injective on $x$ and $x'$.
This reasoning applies to all the points in $\Discaccr$, and the situation of
Proposition~\ref{discaccr_overwhelm_prop} does not occur (every non-injectivity
point is in $\Discaccr$).
\EPR%

\noindent
The finite piecewise increasing case has also the following characteristic: 

\begin{Prop}
  Let $T$ be piecewise increasing on a finite partition. Then $T$ is surjective
  onto~$\II$.
\end{Prop}
\PR
Let $y \in \II$. If $y$ belongs to the images of the continuity
intervals of~$T$ we are trivially done. 
If not, recall that $T(\II)$ has full measure, hence is dense in $\II$, then
similar arguments as in the proof of Proposition~\ref{piecewise_non-inj_prop}
applied to points close to $y$ allow to conclude that $y \in T(\Discaccr)$.
\EPR

Piecewise increasing infinite interval exchanges on a finite partition were
built in~\cite{LN17}\footnote{In~\cite{LN17}, only such interval exchanges were
  studied, explaining why in the definition of infinite interval exchanges
  there, injectivity was imposed on $\II \setminus \Discaccr$.},
and they are the ones whose non-injectivity will be studied in the next
sections.
Nevertheless, note that 
other cases occur.
For instance, let $x \in \Discaccr$ be contained by an interval $[x_1,x_2)$,
  $x_1,x_2 \in \Disc$ for which $T([x_1,x_2))$ is equal either to $[T(x),T(x_2))$
  or $[T(x_1), T(x))$,
so that the graph of $T$ in a neighborhood around $(x,T(x))$ takes respectively
a sprinkled V-like or $\Lambda$-like shape (and $T$ is not piecewise increasing
on a finite partition):

\vspace*{-0.1cm}
\begin{center} 
  \includegraphics[width=5.5cm]{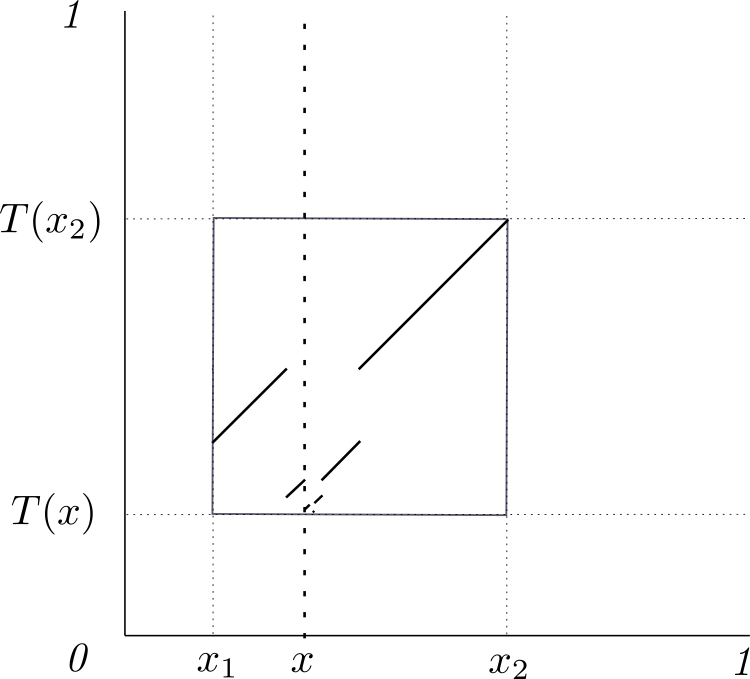}
  \hspace*{1cm}
  \includegraphics[width=5.5cm]{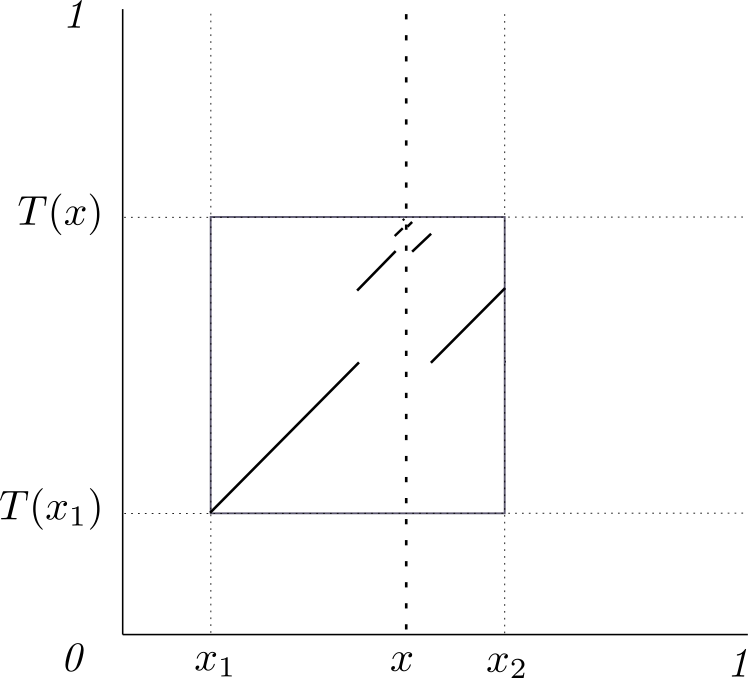}
\end{center}

\vspace*{-0.1cm}
\noindent
Let $J = \II \setminus [x_1,x_2)$. Then, with a configuration in V-like shape in
  $[x_1,x_2)$ as above, $T$ can be such that there is no $x' \in J$ such that
    $T(x')=T(x)$, e.g., when~$T$ has an interval of continuity $[x_3,x_4)
      \subset J$ to the right of $[x_1,x_2)$ such that $T([x_1,x_2) \cup
        [x_3,x_4)) = [T(x_3), T(x_2))$.
Hence, $T$ can be injective while $\Discaccr \neq \emptyset$.
Similarly, with a configuration in $\Lambda$-like shape in $[x_1,x_2)$ as
  above, $T$ can be such that there exists $x' \in J$ such that $T(x')=T(x)$
  with $x' \in \Disc \setminus \Discaccr$, e.g., when~$T$ has an interval of
  continuity $[x_3,x_4) \subset J$ to the right of $[x_1,x_2)$ such that
      $T([x_1,x_2) \cup [x_3,x_4)) = [T(x_1), T(x_4))$.
Hence, $T$ can be non-injective with points of non-injectivity in $\Disc
\setminus \Discaccr$ (only one for each image, in accordance with
Proposition~\ref{discaccr_overwhelm_prop}).
\section{A Construction of Interval Exchange Transformations} 

We recall now how to build families of infinite interval exchange
transformations from symbolic dynamics~\cite{LN17}, and show how the injectivity
of these interval exchanges can be analyzed through this symbolic origin.

\subsection{Basic Definitions for Symbolic Dynamics} 
\label{def_sec}
Let $A$ be a finite alphabet. Let~$A^*$ be the set of finite words over~$A$
including the empty word, and let~$A^\mathbb{N}$ be the set of infinite words
over~$A$ (also called {\smbf sequences}).
If $v \in A^*$ and $w \in A^* \cup A^\mathbb{N}$, their concatenation is denoted
by~$vw$.  A {\smbf factor} of a word~$w$ is a finite word $v$ such that
$w=u'vu''$, where $u', u''$ are possibly empty. Given a word $w$, its set of
factors is denoted by $Fact_w$, and for a set of words $L$, its global set of
factors is $Fact_L=\bigcup_{w \in L} Fact_w$. The subsets of factors of
length~$n$ are respectively denoted by $Fact_w(n)$ and $Fact_L(n)$.
A factor $v$ such that $w=vu$ (resp.~$w=uv$) is a {\smbf prefix} (resp.~{\smbf
  suffix}) of $w$.
We assume that the letters in $A$ are ordered by some order~$\leq$,
lexicographically propagated to all the words in $A^*$ and~$A^\mathbb{N}$.

The set $A^\mathbb{N}$ can be endowed with a topology coming from the {\em
  Cantor metric}: if $w = a_0a_{1}a_{2}...$ and $w' = a'_0a'_1a'_2...$, with
$a_i$, $a'_i \in A$, their distance is~0 if $w=w'$, and~$2^{-k}$ if not, where
$k$ is the smallest non-negative integer for which $a_k \neq a'_k$, i.e., the
length of their larger common prefixes. This topology is equivalent to the
product topology where $A$ is endowed with the discrete topology.
Accordingly, converging sequences in~$A^\mathbb{N}$ are made of words which have
longer and longer common prefixes. One can also embed finite words in $(A \cup
\{\$\})^\mathbb{N}$ by stuffing them to the right with the dummy symbol `$\$$'
making possible to consider converging sequences including~them.
The {\smbf shift map} $\sigma: A^\mathbb{N} \rightarrow A^\mathbb{N}$ is the
erasing of the first letter of its argument, i.e., $\sigma(aw)=w$ for $a\in A, w
\in A^\mathbb{N}$. It is continuous. If $L \subseteq A^\mathbb{N}$ is
topologically closed and invariant under $\sigma$, the pair $(L,\sigma)$ forms a
basic dynamical system in {\em symbolic dynamics}, called a {\smbf subshift} (or
a {\smbf shift space}, or just a {\smbf shift})~\cite{LM95}. Henceforth we
denote it just by $L$.
From a word~$w$ in $A^\mathbb{N}$, a subshift $L$ can be {\smbf induced} as all
the infinite words arising as limits of sequences of words having longer and
longer common prefixes in $Fact_w$, i.e., equivalently, $L =
Closure(\{\sigma^n(w) | n \in \mathbb{N}\})$.
A subshift is {\smbf aperiodic} if it contains no periodic words, i.e., no words
$w$ such that $w=v^\infty$, where $v$ is a finite non-empty word.
If $w\in A^\mathbb{N}$ is {\smbf minimal}, i.e., each factor in~$Fact_w$ occurs
infinitely often in $w$ with bounded gaps, its induced subshift is also {\smbf
  minimal} as a dynamical system, that is, it does not include any proper
non-trivial subshift.  Such a minimal subshift $L$ contains only minimal words
and it can be reobtained from any of them since $Fact_L = Fact_w$ for every $w
\in L$.
Also, a minimal subshift $L$ is {\smbf prolongable},
i.e., each factor in $Fact_L$ can be extended to the right and left so as to
determine factors in $Fact_L$, and also, each word in $L$ can be extended
to the left so as to determine words in~$L$.

A {\smbf substitution} over $A$ is a map $\varphi: A \rightarrow A^*$ which is
extended as a morphism over words as $A^* \cup A^\mathbb{N} \rightarrow A^* \cup
A^\mathbb{N}$ by sending each $w=... a_ia_{i+1}a_{i+2}...$ to
$\varphi(w)=...\varphi(a_i)\varphi(a_{i+1})\varphi(a_{i+2})...$
Accordingly, a substitution can be iterated. 
For instance, let $\varphi$ be defined over $A=\{0,1\}$ as
$\varphi(0) = 01$\label{tm_subst_def} and $\varphi(1) = 10$, then: 
$\varphi(0) = 01$,
$\varphi^2(0) = 0110$, 
$\varphi^3(0) = 01101001$, etc.
\noindent
A {\smbf fixed point} of $\varphi$ is a word $w$ such that
$\varphi(w)=w$.  If there is no letter in $A$ whose image by $\varphi$ is the
empty word, and if $x \in A$ is such that $\varphi(x)=xv$ where $v$ is not the
empty word, then for all $n\geq 0$, $\varphi^{n}(x)$ is a strict prefix of
$\varphi^{(n+1)}(x)$,
so that the sequence $\{\varphi^n(x)\}_{n \in \mathbb{N}}$ converges to the
infinite word $w=\varphi^\infty(x)$ in $A^\mathbb{N}$ which is a fixed point of
$\varphi$.
An {\smbf induced subshift of a substitution}~$\varphi$ is the subshift induced
from one of its fixed points.
In case $\varphi$ is a {\smbf primitive substitution}, i.e., there~exists $n >
0$ such that for every $a_1,a_2 \in A$, $\varphi^n(a_1)$ includes~$a_2$ (see
e.g., the above $\varphi$ example), then all the fixed points of $\varphi$
induce the same subshift $L$, and this subshift is minimal~\cite{Que10}.
\subsection{Special Factors}

Bifurcation possibilities for a subshift $L$ can be described from a classic
word combinatorics-oriented notion~\cite{Cas97,Fog02}:
A factor $v$ in $Fact_L$ is called {\smbf left special} (resp.~{\smbf right
  special}) in~$L$ if $v$ has at least two distinct left (resp.~right) letter
extensions in~$Fact_L$, i.e., $v$ is the suffix (resp.~prefix) of at least two
distinct factors in $Fact_L$.
A {\smbf bispecial factor} is such that it is both right and left special.
Let us consider the following extension sets for~$Fact_L$, $L$ being given:
\[\begin{array}{llll}
Lext(v) &=& \{u \in Fact_L \mid u=av \in Fact_L, \mbox{ with } a \in A\}\\
Rext(v) &=& \{u \in Fact_L \mid u=va \in Fact_L, \mbox{ with } a \in A\}\\
Biext(v) &=&  \{u \in Fact_L \mid u=a_1va_2 \in Fact_L,
                                              \mbox{ with } a_1,a_2 \in A\} 
\end{array}\]
Accordingly, a factor $v$ is respectively left special (resp.~right special)
when $\#Lext(v)>1$ (resp.~$\#Rext(v)>1$), where `$\#$' denotes 
cardinality. 
When $L$ is minimal, every word $w$ in $L$ is such that $Fact_w=Fact_L$, hence
all the above sets can be obtained from a single word $w$ in $L$.

Now, an {\smbf infinite left special word} (or {\smbf infinite left special
  branch})~\cite{Fri05,ACF06,Klo12} of a subshift~$L$ is an infinite word
such that all its prefixes are left special factors in~$Fact_L$.
Since $L$ is a subshift, every such word belongs to $L$, and we denote their set
by $SP_L$.
Let $Lext$ be also defined for $L$ as: 
$$Lext(w) = \{u \in L \mid u=aw \in L, \mbox{ with } a \in A\},$$
The words in $SP_L$ are then the ones such that $\#Lext(w)>1$.  They correspond
to the backward orbit bifurcations in $L$.
If $L$ is an aperiodic subshift induced by a primitive substitution, then $SP_L$
is finite and non empty~\cite{Klo12}. 
Consider its set of left letter extensions\label{lextsp_def}:
$$LextSP_L = \{u \in L \mid u \in Lext(w), \mbox{ with } w \in SP_L\}.$$
Then, this set contains all the words in $L$ on which $\sigma$ is
non-injective.
\subsection{Consecutivity in Subshifts}
\label{consec_sec}

The topology on a subshift $L \subset A^\mathbb{N}$ has a basis made of {\em
  cylinders}, i.e., sets of the form $Cyl_{L}(v) = \{w \in L \mid w=vu\}$ for
each $v \in Fact_{L}$, which are {\em clopen sets}. Subshifts are indeed
zero-dimensional spaces.

Recall then that we assume the alphabet $A$ of $L$ as being ordered
by~`$\leq$', lexicographically extended to all the words over $A$.
Two words $w$ and $w'$ in a subshift $L$ are then said to be {\smbf consecutive}
if $w<w'$ and if there is no word $w'' \in L$ such that $w < w'' < w'$.
Intervals of words in $L$ are defined as $[w, w'] = \{w'' \in L \mid w \leq w''
\leq w'\}$.
Induced by~`$\leq$' there is a strict partial order over the intervals of $L$
where $[w,w'] < [u,u']$ if $w' < u$.
So, two intervals $[w,w']$, $[u,u']$ are {\smbf consecutive} if
$w'$ and $u$ are consecutive. 
Every $Cyl_{L}(v)$ is an interval of words in~$L$, and $Cyl_{L}(v) <
Cyl_{L}(v')$ if $v < v'$ with $v' \notin Rext(v)$.
By compactness of $L$ as a subshift, each~$Cyl_{L}(v)$ has two endpoints in~$L$
as its smallest and its greatest words. Hence, $Cyl_{L}(v)$ and $Cyl_{L}(v')$ are
consecutive if the greatest word in $Cyl_{L}(v)$ is consecutive to the smallest
word in $Cyl_{L}(v')$.
Note that for a given $n>0$, the set of cylinders corresponding to the set
$Fact_L(n)$, i.e., the factors of length~$n$ in $L$, induces an ordered
partition of~$L$ as $\bigsqcup_{v \in Fact_L(n)} Cyl_{L}(v)$ of consecutive
cylinders.

We denote by $\ShiftsAM$ the set of all the aperiodic minimal subshifts. 
Here are a few properties about consecutivity in these subshifts:
\begin{Lem}\label{conseq_pairs_only_lem}
  Let $L \in \ShiftsAM$.  Consecutive words in $L$ occur only in pairs.
\end{Lem}
\PR $L$ being aperiodic minimal, it has all its cylinders made of an infinite
number of words. Hence, there cannot be $w<w''<w'$, all consecutive.
\EPR

\noindent
In view of the above lemma we define $\ConsecPairL \subset L \times L$ as the
set of pairs $(w,w')$ such that $w < w'$ are consecutive words in $L$, and
$\ConsecL \subset L$ as the set of words that belong to a pair in
$\ConsecPairL$.
The full link between consecutive words and cylinders is then~\cite{LN17}:
\begin{Lem}\label{conseq_small_great_lem}
  Let $L \in \ShiftsAM$. A pair  $(w, w')$ belongs to $\ConsecPairL$ iff
  $w$ and $w'$ are the greatest and the smallest words of two consecutive
  cylinders.
\end{Lem}

\noindent
From the above lemma, we readily have:

\begin{Lem}\label{double_conv_seq_lem} 
  Let $L \in \ShiftsAM$.  Then $w \notin \ConsecL$ iff $w$ is the limit of
  both increasing and decreasing sequences of words in~$L$.
\end{Lem}

\noindent
Now, with regard to the shift $\sigma$, non-consecutivity can be preserved: 

\begin{Lem}\label{non-consec_preserv_lem} 
  Let $L \in \ShiftsAM$.
  Let $aw \in L$, with $a \in A$. If $aw \notin \ConsecL$, then $w \notin
  \ConsecL$.
\end{Lem}
\PR According to Lemma~\ref{double_conv_seq_lem}, $aw$ is the limit of both
increasing and decreasing sequences of words in~$L$. An increasing
(resp.~decreasing) sequence $\{aw^{(i)}\}$ of words converging to $aw$ is such
that the $aw^{(i)}$ have longer and longer common prefixes as $i$ goes to
infinity.
This property remains true for $\{\sigma(aw^{(i)})\}=\{w^{(i)}\}$ which
converges to $w$, and these sequences remain increasing (resp.~decreasing) since
$\sigma$ is order-preserving in $Cyl_L(a)$.
\EPR

\noindent
Note that the converse does not hold: even with a prolongable $L$, if $w \notin
\ConsecL$, all the sequences $\{w^{(i)}\}$ converging to $w$ cannot necessarily
be turned into sequences $\{aw^{(i)}\}$ converging to $aw$, 
since the extension letters of $w^{(i)}$ can be different from~$a$.

About the non-consecutivity preservation by $\sigma$, it is also sometimes
convenient to consider the more constrained set $\ConsecPairL$:

\begin{Lem}\label{non-consec_pair_preserv_lem} 
  Let $L \in \ShiftsAM$.
  Let $aw$, $aw' \in L$, with $a \in A$. If $(aw, aw') \notin \ConsecPairL$,
  then $(w, w') \notin \ConsecPairL$.
\end{Lem}
\PR If $(aw, aw') \notin \ConsecPairL$, there is $aw''\in L$ such that
$aw < aw'' < aw'$. Hence, $w < w'' < w'$ since $\sigma$ is order-preserving in
$Cyl_L(a)$.
\EPR

\noindent
Again the converse does not hold: if there is $w'' \in L$ such that $w < w'' <
w'$, then $w''$ does not necessarily extend as $a w''$ so that $aw < aw'' <
aw'$.
\subsection{Going from Subshifts to $\II$} 
\label{T_L_construct_sec}

We now summarize how subshifts $L$ can be semi-conjugated to functions
over $\II$ as {\em interval exchange transformations} (generally
infinite)~\cite{LN17}.

First of all, when $L \in \ShiftsAM$ it can be endowed with a $\sigma$-invariant
and non-atomic Borel probability measure~$\mu$. This measure is positive on the
cylinders $Cyl_{L}(v)$ for each $v \in Fact_{L}$, and these generate the
corresponding Borel $\sigma$-algebra~\cite{LM95}.  Since $\mu$ is non-atomic it
is null on single words.
Such a measure can be based on
the occurrence frequencies of the factors in~$L$, and when $L$ has been induced
by a primitive substitution, explicit formulas exist to compute these
frequencies~\cite{Que10}.
Now, let $L$ be an ordered subshift in $\ShiftsAM$ measured by $\mu$, and
consider the following valuation map:

\[\begin{array}{lllll} 
 \phi_\mu: \;&L \;\;\;\; &\rightarrow \;\;\;\; & [0,1]\\
             &  w & \mapsto &\mu([w_{min}, w])
\end{array}\]

\noindent
where $w_{min}$ is the smallest word of $L$ (which exists since $L$ is
closed and bounded). The map $\phi_\mu$ can be proved to be monotonic
non-decreasing, continuous, surjective and {\em non-injective on pairs
  of consecutive words}.

The map $\phi_\mu$ is then made injective on $\LprimOne \subset L$, 
defined such that
if $\phi_\mu^{-1}(x)$ is made of two consecutive words in $L$ (see
Lemma~\ref{conseq_pairs_only_lem}), only the greater one belongs to
$\LprimOne$, that is, 
we get a well-defined induced inverse map $\phiinv$ as:
\[\begin{array}{lllll}\label{phiinv_def} 
 \phiinv: \;&[0,1] \;\;\;\; &\rightarrow \;\;\;\; & \LprimOne\\
             &  x & \mapsto &\left\{
              \begin{array}{lll}
                \phi_\mu^{-1}(x) & \mbox{\em if $\phi_\mu(x)$ is invertible}\\
                w' & \mbox{\em if $\phi_\mu^{-1}(x)=\{w,w'\}$, $w < w'$.}
             \end{array} \right.
\end{array}\]
Since consecutive words $\{w,w'\}$ only occur as endpoints of cylinders (see
Lemma~\ref{conseq_small_great_lem}), $\LprimOne$ has full measure in $L$.
Note also that to consider partitions of $L$ leading by $\phi_\mu$ to partitions
of $[0, 1]$, one must consider $\LprimOne$ instead of $L$ since otherwise
consecutivity would lead to overlaps in $[0, 1]$.

We can then define the following map:
\[\begin{array}{lllll} 
 T_L: \;&\II \;\;\;\; &\rightarrow \;\;\;\; & \II\\
        &  x & \mapsto &\phi_\mu(\sigma(\phiinv(x))) 
\end{array}\]
whose domain is set to $\II$ rather than $[0,1]$ since $T_L$ can be proved 
right-continuous, and $T_L(\II)\subset \II$.
Accordingly we put $\Lprim = \LprimOne \setminus \{\phiinv(1)\}$.
The map $T_L$ can be then shown to preserve Lebesgue measure, and when the
measure $\mu$ of the set of the infinite left special words in $L$ is null,
i.e., $\mu(SP_L) = 0$,
then $T_L$ can be proved to be an {\em interval exchange transformation}.
Moreover, $\phi_\mu$ is a conjugacy between $T_L$ and $\sigma$ on $\Lprim$, that
is:
$$ \phi_\mu(\sigma(w))= T_L(\phi_\mu(w)), \;\; \forall w \in \Lprim,$$
which becomes a topological semi-conjugacy when considered on all~$L$. 

The construction of $T_L$ can be made explicit by building a partition of $\II$
for which $T_L$ is a local isometry on the components of positive measure.
The idea is to start from the symbolic side and find word intervals in the
subshift $L$ on which $\sigma$ is injective.
Indeed, by measure preservation we always have 
$\mu(\sigma^{-1}(\sigma([w,w']))) =  \mu(\sigma([w,w']))$.
Moreover, when $\sigma$ is injective on $[w,w']$ then
$\sigma^{-1}(\sigma([w,w'])) = [w,w']$ and $\sigma([w,w']) =
  [\sigma(w),\sigma(w')]$.
Hence, for every word $w''$ in $[w,w']$, we have
$\phi_\mu(\sigma(w'')) - \phi_\mu(w'') = (\mu([w_{min},\sigma(w)] +
\mu([\sigma(w),\sigma(w'')])) - (\mu([w_{min},w]) + \mu([w,w''])) =
\mu([w_{min},\sigma(w)]) - \mu([w_{min},w])$, 
that is, on $\phi_\mu([w,w'])$ the map $T_L$ is a translation which is the
expected property for a continuity interval of $T_L$.

Now, when $L$ is minimal it is prolongable, and a criterion for~$\sigma$ to be
invertible on an interval $[w,w']$ is that $w$ and $w'$ have a common prefix
$xu$ such that $x \in A$ and $u$ is not a left special factor.
Thus, $\sigma$ is injective over each cylinder component of a partition of $L$
of the following form\label{part_L_def}:
$$  PART_{L}: \;\;\;\; L =
    \bigsqcup_{k >0} Cyl_{L}(v^{(k)}) \;\; \sqcup \;\; LextSP_{L},$$
where each $v^{(k)}=x^{(k)}u^{(k)}$ is such that $x^{(k)} \in A$ and $u^{(k)}$
is not a left special factor, and where $LextSP_{L}$ is the set of all the left
extensions of the infinite left special words in $SP_{L}$ (see
p.~\pageref{lextsp_def}), which accordingly cannot belong to any cylinder of the
form $Cyl_{L}(v^{(k)})$.
In order to build such a partition $PART_L$ one can use the fact
that a non left special factor~$v$ has all its right extensions as non left
special too, and apply the following iterative process:

\begin{enumerate}[i.]\itemsep=0cm\label{part_L_algo}
\item For each factor $a_{i_1}a_{i_2}$ in $Fact_L(2)$, with $a_{i_1},a_{i_2} \in
  A$, if $a_{i_2}$ is not left special, select this factor, and add
  $Cyl_{L}(a_{i_1}a_{i_2})$ to $PART_{L}$.

\item For each factor $a_{i_1}a_{i_2}$ not selected in the preceding step,
   consider each factor $a_{i_1}a_{i_2}a_{i_3}$ in $Fact_L(3)$, 
   and if $a_{i_2}a_{i_3}$ is not left special, select this factor, and add
   $Cyl_{L}(a_{i_1}a_{i_2}a_{i_3})$ to $PART_{L}$.

\item Do the same as above for every length $n$ and for each factor
  $a_{i_1}a_{i_2} \cdots a_{i_{n-1}}$ which was not selected in the preceding
  $(n-1)$th step, that is, such that $a_{i_2} \cdots a_{i_{n-1}}$ is still left
  special.

\item After having carried on the preceding iterations for every $n$, the
  remaining words can only be the ones making $LextSP_{L}$.
\end{enumerate}

\noindent
From $PART_{L}$ one can then build a partition $PART_{L,\II}$ of $\II$: Let us
first define the map which gives the image in $\II$ by $\phi_\mu$ of any subset
in $L$ but using $\Lprim$ only:

\[\begin{array}{lllll} 
\phiLprim: \;& 2^L \;\;\;\; &\rightarrow \;\;\;\; & 2^{\II}\\
       &  Y \subset L & \mapsto & \{x \in \II \; | \; \phiinv(x) \in Y\}.
\end{array}\]

\noindent
Accordingly, $\phiLprim(Y) = \phiLprim(Y \cap \Lprim)$, and $\phiLprim(\Lprim) =
\II$. Also, if $w \notin \Lprim$, $\phiLprim(\{w\})= \emptyset$.
Next, we just apply $\phiLprim$ to the components of $PART_L$:
$$  PART_{L,\II}: \;\;\;\; \II =
    \bigsqcup_{k > 0} \phiLprim(Cyl_{L}(v^{(k)})) \;\; \sqcup \;\;
    \phiLprim(LextSP_{L}).$$
By construction of $PART_L$ and its components $Cyl_{L}(v^{(k)})$, the map $T_L$
is an order-preserving isometry on each interval $\phiLprim(Cyl_{L}(v^{(k)}))$,
and the accumulation points of its discontinuities can only occur in
$\phiLprim(LextSP_{L})$.
Moreover, according to Section~\ref{iiets_sec}, we know that these points are
where non-injectivity of $T_L$ takes its source.
We precise this fact further in the next section.

\subsection{Non-Injectivity of the Built Interval Exchange Transformations}
\label{non-inj_iets_sec}

From now on we shall only deal with ordered minimal subshifts $L$ in $\ShiftsAM$
equipped with a $\sigma$-invariant measure $\mu$ positive on the cylinders such
that $\mu(SP_L) = 0$, so that $T_L$ is an interval exchange transformation, and
we denote their set by $\ShiftsAMZ$.
We also denote by $\Ninj \subset \II$ the set where $T_L$ is non-injective.
Then, precising the intuition expressed at the end of the last section, and
considering the set $\phiLprim(LextSP_{L})$, we have that $T_L$'s
non-injectivity on $\II$ is indeed tied to $\sigma$'s non-injectivity on~$L$:
\begin{Lem}\label{ninj_inclusion_lem}
  Let $L \in \ShiftsAMZ$,
  Then
  $\Ninj \subset \phiLprim(LextSP_L).$
\end{Lem}
\PR Since $T_L(x) = \phi_\mu(\sigma(\phiinv(x)))$, and since $\phiinv$ is
injective, the only way for $T_L(x)$ to be non-injective is for points in
$\phiLprim(LextSP_L)$ because $LextSP_L$ is the set of words in $L$ for which
$\sigma$ is non-injective.
\EPR

\noindent
Note that since~$T_L$ is built from~$\Lprim$ only, strict inclusion may hold in
the above lemma. Here is a characterization of the set $\Ninj$ considering
non-injectivity on pairs of its points:

\begin{Lem}\label{noninj_lem}
  Let $L \in \ShiftsAMZ$.
  Then $x, x' \in \Ninj$ with $T_L(x) = T_L(x')$
  iff there is $w \in SP_L$ with $\phiLprim(\{a_1w, a_2w\}) =
  \{x,x'\}$, $a_1,a_2 \in A$, $a_1w, a_2w \in \Lprim$.
\end{Lem}
\PR By definition of $\Lprim$, and according to
Lemma~\ref{conseq_pairs_only_lem},
if $a_1w$, $a_2w \in \Lprim$, they cannot be consecutive.
Hence there exists a cylinder $C$ such that $a_1w < C < a_2w$, and
since $\mu(C)>0$ its image $\phi_\mu(C)$ has positive measure too,
so that $x=\phi_\mu(a_1w) < x'=\phi_\mu(a_2w)$.
But then $T_L$ is non-injective on them since
$T_L(x)=T_L(\phi_\mu(a_1w)) = \phi_\mu(\sigma(a_1w)) = \phi_\mu(w)
=\phi_\mu(\sigma(a_2w)) = T_L(\phi_\mu(a_2w)) = T_L(x').$
Conversely, if $x, x' \in \Ninj$ with $T_L(x) = T_L(x')$, then by
Lemma~\ref{ninj_inclusion_lem}, there exist $a_1w, a_2w \in LextSP_L$ with
$\phiLprim(a_1w)=x$ and $\phiLprim(a_2w)= x'$. By definition of $\phiLprim$ both
$a_1w$ and $a_2w$ must belong to $\Lprim$.%
\EPR

\noindent
A main case where the conditions of Lemma~\ref{noninj_lem} are not satisfied is
when $\phiLprim(\{a_1w, a_2w\}) = \{x\}$ meaning that $(a_1w,a_2w) \in
\ConsecPairL$, hence that only $a_2w$ belongs to $\Lprim$ and gives~$x$: No
non-injectivity is then induced from $a_1w$ and $a_2w$.
Moreover if $Lext(w)$ is here exactly $\{a_1w, a_2w\}$, i.e., $\#Lext(w)=2$,
then $x$ is ensured to not belong to $\Ninj$ since there is no other word $w'$
in $LextSP_L$, hence in $L$, for which $\phi_\mu(\sigma(w')) = x$.  This is a
situation where $\Ninj \neq \phiLprim(LextSP_L)$ in
Lemma~\ref{ninj_inclusion_lem}.

Note that to be in position to apply Lemma~\ref{noninj_lem} one can use
Lemma~\ref{double_conv_seq_lem} or~\ref{non-consec_preserv_lem} about
non-consecutivity, since if $w \notin \ConsecL$ then $w \in \Lprim$ in
accordance with the definition of $\phiinv$ (see p.~\pageref{phiinv_def}).
Also, to establish that two words $w$, $w'$ are not both in $\Lprim$, one can
apply Lemma~\ref{conseq_small_great_lem} to prove that they are consecutive.
Moreover, Lemma~\ref{non-consec_pair_preserv_lem} may be sometimes used to prove
consecutivity in $LextSP_L$ through consecutivity in~$SP_L$.

Also, with regard to the structure of an interval exchange $T_L$ built from a
subshift $L$,
it is worth to note an additional property about $\Ninj$:

\begin{Lem}\label{ninj_discaccr_lem}
  Let $L \in \ShiftsAMZ$.
  Then $\Ninj = \Discaccr$, where $T_L$ is piecewise increasing on a finite
  partition.
\end{Lem}
\PR Consider the partition $L =\bigsqcup_{a \in A} Cyl_{L}(a)$.
Then $T_L$ is increasing on each $\phiLprim(Cyl_L(a))$, $a \in A$, since
$\sigma$ is order-preserving on $Cyl_{L}(a)$ and since $\phi_\mu$ is strictly
increasing on $L$.
Hence, $T_L$ is piecewise increasing on a finite partition, and by
Proposition~\ref{piecewise_non-inj_prop}, $\Ninj = \Discaccr$.
\EPR

Recall that the construction of $PART_L$ in Section~\ref{T_L_construct_sec}
yields converging sequences to words in $LextSP_L$, made of endpoints of
cylinders corresponding to continuity intervals in $T_L$.
Accordingly, these endpoints are sent by $\phi_\mu$ to interval endpoints of the
partition $PART_{L,\II}$ that accumulate onto points belonging to
$\phiLprim(LextSP_L)$.
However, we know from Lemma~\ref{noninj_lem} that $\phiLprim(LextSP_L)$ is not
necessarily equal to $\Ninj$, that is, $PART_{L,\II}$ is not necessarily maximal
with respect to the discontinuities of $T_L$. Hence note that $PART_{L,\II}$
cannot straightforwardly be used to ensure that a point belongs to~$\Discaccr$,
and thus by Lemma~\ref{ninj_discaccr_lem}, to~$\Ninj$.
\section{A Generalized Thue-Morse Case}

In this section, from a particular set of subshifts, we explicitly build a
family of infinite interval exchange transformations $T_L$ having a
characterizable finite non-injectivity set $\Ninj$.

\subsection{Thue-Morse Sequences}
\label{thue-morse_sec}

The {\smbf Thue-Morse sequence} $t=t^{(0)}t^{(1)}t^{(2)}...$, where the $t^{(i)}$
are letters in $A=\{0,1\}$, is a classic infinite word over two letters with a
long history~\cite{Pro1851,Thu12,Mor21b}
with several ways of being defined~\cite{AS98,AS03}:
Its number theoretic-oriented definition is $t^{(i)} = s_2(i) \;mod\; 2$, where
$s_2$ is the base-2 sum-of-digits function which counts the number of $1$'s in
the base-$2$ representation of~$i$.
Its recursive definition is $t^{(0)}=0$, $t^{(2n)}=t^{(n)}$ and
$t^{(2n+1)}=1-t^{(n)}$, $n \geq 0$ (that is, $1-t^{(n)}$ is the {\em mirror} word
of~$t^{(n)}$).
A third definition relies on a substitution (the same as in the example
p.~\pageref{tm_subst_def}), i.e., 
\[\begin{array}{llllll}
  \varphi(0) = 01, \;\;\;  \varphi(1) = 10, 
\end{array}\]  
for which $t$ is obtained as the fixed point $\varphi^\infty(0)$.

The Thue-Morse sequence $t$ has been generalized in several ways, in particular
for greater alphabets.  Here, we focus on the one given as a variation of the
number theoretic-oriented definition by taking a base-2 sum-of-digits function
over $\mathbb{Z}_m$ with $m \geq 2$, instead of $m=2$, that is, $t^{(i)}_m =
s_2(i) \;mod\; m$,
determining the {\smbf generalized Thue-Morse sequence} $t_m$ over $m$ letters
$A_m=\{0,1,...,(m-1)\} = \mathbb{Z}_m$ \cite{TS95,AS03}
(considering the alphabet as~$\mathbb{Z}_m$ is often convenient since letters
can be then expressed using the structure of~$\mathbb{Z}_m$).
These sequences also have an equivalent substitution-oriented definition
where~$t_m$ is obtained as the fixed point $\varphi_m^\infty(0)$ where:
\[\begin{array}{llllll}
   \varphi_m(k) = k (k+1),  & \mbox{for every } k \in A_m=\mathbb{Z}_m.
\end{array}\]  

\noindent
For each $m$, the generalized Thue-Morse sequence $t_m$ induces what we call the
{\smbf Thue-Morse subshift}~$L_m$, lexicographically ordered from the natural
order over $\mathbb{Z}_m$. Since $\varphi_m$ is primitive, these subshifts are
minimal~\cite{Que10}.
Moreover, $t_m$ is known to be aperiodic for every~$m$~\cite{MM91,AS00},
hence $L_m \in \ShiftsAM$. 
Note also that $\varphi_m$ has exactly $m$ distinct fixed
points\label{fix_varphi}:
$$ Fix(\varphi_m) = \{w \in A_m^\mathbb{N} \mid
                 w=\varphi_m^\infty(k),\; k \in A_m\}.$$ 
By minimality, each word in $Fix(\varphi_m)$ has the same set of factors, so
that each induces~$L_m$, not only $\varphi_m^\infty(0)$.
Also, $\varphi_m$ being primitive, each $L_m$ is {\em uniquely
  ergodic}~\cite{Que10},\label{uniq_ergod}
and the corresponding $\sigma$-invariant probability measure~$\mu$ can be
explicitly obtained from its factor frequencies as computed in~\cite{Bal12}.
\subsection{The Special Words of $L_m$}
\label{special_words_sec}

Following Section~\ref{non-inj_iets_sec},
to study the non-injectivity of the conjugated interval exchanges $T_{L_m}$, we
need to describe the left letter extension set $LextSP_{L_m}$, and thus
beforehand the infinite left special words of $SP_{L_m}$.

First of all, the finite special factors of the Thue-Morse sequences $t_m$,
hence equivalently the ones of their associated subshifts $L_m$, have been
already studied.
In particular the next lemma can be deduced from~\cite{Sta12}:

\begin{Lem}\label{bisp_words_lem}
  Let $L_m$ be a Thue-Morse subshift, $m\geq 2$. 
 Then its bispecial factors have the following possible three forms and extension
 cardinalities:
 \begin{itemize}\itemsep=0.1cm
 \item The single letters $k$ in $A_m=\mathbb{Z}_m$, for which $\#Biext(.) = 2m-1$.
 \item The length-2 words $k(k+1)$, with $k \in A_m=\mathbb{Z}_m$, and all their
   iterates $\varphi_m^n(k (k+1))$, $n\geq 0$, for which $\#Biext(.) = 2m$.

 \item The length-3 words $k(k+1)(k+2)$, with  $k \in A_m=\mathbb{Z}_m$, and all their
   iterates $\varphi_m^n(k (k+1)(k+2))$, $n\geq 0$, for which
   $\#Biext(.) = 2$.
 \end{itemize}
\end{Lem}
\PR According to~\cite{Sta12}[Lemma 4.1], every bispecial factor $v \in
Fact_{L_m}$ of length $\geq 4$ is such that there is another bispecial factor $u
\in Fact_{L_m}$ with $v=\varphi_m(u)$, also preserving cardinality, i.e.,
$\#Biext(v) = \#Biext(u)$. Next, the initial bispecial factors of length $\leq
3$ are shown in~\cite{Sta12}[proof of Lemma 4.2] to be $k$, $k(k+1)$ and
$k(k+1)(k+2)$, $\forall k \in \mathbb{Z}_m$, having the above indicated
cardinalities.
\EPR

\noindent
The ``{\em bispecial legacy property}'' implicitly on use in the above proof is
that a bispecial factor in $L_m$ remains bispecial under $\varphi_m$.  Indeed,
since all the first (resp.~last) letters of the images $\{\varphi_m(k)\}_{k\in
  A_m}$ are distinct, if $v$ is bispecial, every distinct extension $a_1va_2$,
with $a_1,a_2 \in A$, is such that
$\varphi_m(a_1va_2)=\varphi_m(a_1)\varphi_m(v)\varphi_m(a_2)$ yields a distinct
extension for $\varphi_m(v)$ (this is just the simplest situation to obtain most
of the bispecial factors of a subshift induced by a substitution~\cite{Klo12}).
From this remark some of the associated cardinalities of these bispecial factors
can also be precised:

\begin{Lem}\label{sp_degree_lem}
Let $L_m$ be a Thue-Morse subshift, $m\geq 2$. 
Every bispecial factor $v \in Fact_{L_m}$ is such that:
 \begin{itemize}\itemsep=0.1cm
   \item When $\#Biext(v) = 2$, then $\#Lext(v)=\#Rext(v)=2$.
   \item When $\#Biext(v) = 2m-1$ or $2m$, then $\#Lext(v)= \#Rext(v)=m$.
 \end{itemize}
\end{Lem}
\PR If $\#Biext(v)=2$, then $\#Lext(v)$ and $\#Rext(v)$ cannot be $>2$ since
otherwise $\#Biext(v)$ would also be $>2$, $L_m$ being prolongable.  This
cardinality cannot be $<2$ since $v$ would trivially not be bispecial.

For the other cases, note first that by Lemma~\ref{bisp_words_lem}, $\forall k
\in \mathbb{Z}_m$, $k(k+1) \in Fact_{L_m}(2)$.  But we also have that $\forall
k,i \in \mathbb{Z}_m$, $k(k+i) \in Fact_{L_m}(2)$ since $\varphi_m^h(k(k+1))$ is
also in $Fact_{L_m}$ with $\varphi_m^h(k(k+1)) = \varphi_m^h(k)
\varphi_m^h(k+1)$ where $\varphi_m^h(k)$ ends with $k+h$, and $\varphi_m^h(k+1)$
starts with $k+1$, while $k$ and $h$ can vary independently.
Hence, $\#Lext(k)= \#Rext(k)=m$, $\forall k \in \mathbb{Z}_m$.
Then, the {\em bispecial legacy property} on $\{\varphi_m^n(k)\}_{k \in
  \mathbb{Z}_m}$ for $n \geq 0$ gives the result.\EPR

\begin{Lem}\label{sp_fix_lem}
  Let $L_m$ be a Thue-Morse subshift, $m\geq 2$. 
  Then $$SP_{L_m} = Fix(\varphi_m),$$
  and for every $w \in SP_{L_m}$, we have $\#Lext(w)=m$. 
\end{Lem}
\PR Every fixed point $w = \varphi_m^\infty(k)$ belongs to $SP_{L_m}$ since,
according to Lemma~\ref{bisp_words_lem}, $\varphi_m^n(k)$ is a bispecial factor
for every $n\geq 0$, and $\varphi_m^{n+1}(k)=
\varphi_m^n(k)\varphi_m^n(k+1)$. Hence, $\varphi_m^n(k)$ is a strict
left-special prefix of $\varphi_m^{n+1}(k)$, thus also of $w$,
and prefixes of a left special prefix are left special too. Moreover, according
to Lemma~\ref{sp_degree_lem}, all these prefixes are such that $\#Lext=m$,
hence $\#Lext(w)=m$ since $L_m$ is closed and one can build converging
distinct sequences to each word in $Lext(w)$.

Conversely, if $w \in SP_{L_m}$, all its prefixes are left special, and there
are also arbitrarily long factors among these prefixes which must be right
special, that is, bispecial.  Indeed, if there would be a length up to which no
prefix of $w$ is right special in $L_m$, then $w$ would be univocally
determined, hence ultimately periodic, contradicting that $L_m \in \ShiftsAM$.
Now, there are at least two left letter extensions $a_1w$ and $a_2w$ of $w$.
Each must also have infinitely many prefixes which are right special, but
then also their suffixes which are prefixes of~$w$.
This makes infinitely many bispecial prefixes of~$w$ with $\#Biext > 2$. By
Lemma~\ref{bisp_words_lem}, these must have the form $\varphi_m^n(k)$ where $k$
is the first letter of $w$ since they are the only ones starting by $k$.  Since
$\varphi_m^n(k) = \varphi_m^{n-1}(k)\varphi_m^{n-1}(k+1)$, $\forall n>0$, then
$w=\varphi_m^\infty(k)$. 
\EPR

\noindent
As we have noted in p.~\pageref{fix_varphi}, we have $\#Fix(\varphi_m)=m$ so
that: 
\begin{Cor}
  Let $L_m$ be a Thue-Morse subshift, $m\geq 2$. Then
  $$\mu(SP_{L_m})=0.$$
\end{Cor}

\noindent
Thus, the subshifts $L_m$ are in $\ShiftsAMZ$ for every~$m$, and all fit to be
associated with a semi-conjugated~$T_{L_m}$ that can be called a {\smbf
  Thue-Morse interval exchange transformation}.
\subsection{Non-Injectivity of $T_{L_m}$} 

From Lemma~\ref{sp_fix_lem}, since $\#Lext(w)=m=\#A$ for all $w \in SP_{L_m}$,
the set of the left letter extensions of the infinite left special words is:
$$LextSP_{L_m} = \{ w' \in L_m \mid w' = a w, \mbox{ with } a \in A_m, w \in
SP_{L_m}\}.$$

\noindent
A first property about this set is an extension and converse of
Lemma~\ref{non-consec_pair_preserv_lem}:
\begin{Lem}\label{consec_preserv_TM_lem} 
  Let $L_m$ be a Thue-Morse subshift, $m\geq 2$.  Let $w$, $w' \in SP_{L_m}$.
  Then $(w, w') \in \ConsecPairLm$ iff $(aw, aw') \in \ConsecPairLm$,
  $\forall a \in A_m$. 
\end{Lem}
\PR
Every $w \in SP_{L_m}$ is such that $aw \in L_m$ for each $a \in A_m$, and
Lemma~\ref{non-consec_pair_preserv_lem} applies to any pair of them.
Conversely, if $w,w' \in SP_{L_m}$ with some $w''\in L_m$ such that $w<w''<w'$,
since $L_m$ is prolongable, there is some $a \in A_m$ such that $aw'' \in L_m$,
so that $aw<aw''<aw'$.
\EPR

\noindent
Also for Thue-Morse subshifts $L_m$, a converse of
Lemma~\ref{non-consec_preserv_lem} holds for $SP_{L_m}$:

\begin{Lem}\label{non-consec_preserv_TM_lem} 
 Let $L_m$ be a Thue-Morse subshift, $m\geq 2$, 
 and let $w \in SP_{L_m}$.
 If~$w \notin \ConsecLm$, then $aw \notin \ConsecLm$,
  $\forall a \in A_m$. 
\end{Lem}
\PR Following Lemma~\ref{double_conv_seq_lem} we prove that if there exist both
increasing and decreasing sequences of words converging to $w$, this is also the
case for all its extensions $aw$, $a \in A_m$.
So let ${\cal S} = \{v_n\}_{n \in \mathbb{N}}$ be an increasing
(resp.~decreasing) sequence of words in $L_m$ converging to $w$. Since $L_m$ is
prolongable, up to extracting a subsequence from ${\cal S}$ there is at least
one letter $i \in A_m$ such that ${\cal S}_i=\{iv_n\}_{n \in \mathbb{N}}$
is also an increasing (resp.~decreasing) sequence of words in $L_m$ converging
to $iw$.
By the form of $\varphi_m$, we have $\varphi_m(iv_n) =
\varphi_m(i)\varphi_m(v_n) = i(i+1)\varphi_m(v_n)$.
This means that from ${\cal S}_i$, we can obtain another sequence in $L_m$ as
${\cal S}_{i+1}= \{(i+1)\varphi_m(v_n)\}_{n \in \mathbb{N}}$ converging to
$(i+1)\varphi_m(w)$, i.e., $(i+1)w$, since $w$ is in $SP_{L_m}$ and is a fixed
point of $\varphi_m$ according to~Lemma~\ref{sp_fix_lem}.
Now, $\varphi_m$ is an {\em increasing map}, i.e., for every $w, w' \in L_m$, if
$w < w'$ then $\varphi_m(w) < \varphi_m(w')$.
Thus, the sequence ${\cal S}_i$ being increasing (resp.~decreasing) means that
${\cal S}_{i+1}$ is increasing (resp.~decreasing) too.
Applying $\varphi_m^j({\cal S}_i)$ for each $0 \leq j < m$ yields sequences
${\cal S}_{i+j}$ converging to $(i+j)w$, that is, to every letter extension
making $Lext(w)$. Hence, if $w \notin \ConsecLm$, every word in $Lext(w)$ is not
in $\ConsecLm$ either.  \EPR

Thus, consecutivity among the words in $LextSP_{L_m}$ of the Thue-Morse
subshifts $L_m$ --~hence their non-injectivity~--, can be studied from
consecutivity among the words in~$SP_{L_m}$, which is much easier because of
the results in Section~\ref{special_words_sec}.
From now on, let us denote the $m$ words of $SP_{L_m}$ as $w_0,\cdots,w_{m-1}$,
where $w_i = \varphi_m^{\infty}(i)$ for each $i \in A_m$ (see
Lemma~\ref{sp_fix_lem}):

\begin{Lem}\label{conseq_sp_lem}
  Let $L_m$ be a Thue-Morse subshift, $m\geq 2$. 
  In $SP_{L_m}$, we have $(w_{m-2}, w_{m-1}) \in \ConsecPairLm$, and 
  when $m>2$, $w_0,..., w_{m-3} \notin \ConsecLm$.
\end{Lem}
\PR Following Lemma~\ref{conseq_small_great_lem}, let us first show that
$w_{m-2}$ is the greatest word in $Cyl_{L_m}(m-2)$, and $w_{m-1}$ is the
smallest word in $Cyl_{L_m}(m-1)$ (recall that $m-2$ and $m-1$ are the two last
letters of $A_m = \mathbb{Z}_m$).
About $w_{m-2}$, assume on the contrary that there would exist a word $w \in
Cyl_{L_m}(m-2)$, such that $w_{m-2} < w$. Let $v$ be their maximal common prefix
so that $w_{m-2}=vu$ and $w=vu'$, with $u,u' \in L_m$.
Then, $v$ must be bispecial since every prefix of a word in $SP_{L_m}$ is
left-special, and $u$ and $u'$ must start with different letters, respectively
$a$ and $a' \in A_m$.
Since $va$ is left special being a prefix of $w_{m-2}$, and since $va'$ is
extendable to the left ($L_m$ being prolongable),
we have that $\#Biext(v) > 2$.
Now, according to Lemma~\ref{bisp_words_lem}, being a prefix of $w_{m-2}$, that
is $\varphi_m^{\infty}(m-2)$, the factor~$v$ must necessarily be of the form
$\varphi_m^n(m-2)$ for some $n \geq 0$.
Then by its construction, $w_{m-2}$ is such that $\varphi_m^{n+1}(m-2)$ is one
of its prefixes too.
Since $\varphi_m^{n+1}(m-2)=\varphi_m^n((m-2) (m-1)) =
\varphi_m^n(m-2)\varphi_m^n(m-1)$,
the extension letter $a$ to the right of $\varphi_m^n(m-2)$ must be $m-1$, i.e.,
$v (m-1)$ is a prefix of $w_{m-2}$. But then there is no letter $a'$ in $A_m$
which is greater than $m-1$.  Hence, there is no word greater than $w_{m-2}$
such as $w$ in $L_m$.

Similar arguments apply to show that $w_{m-1}$ is the smallest word in
$Cyl_{L_m}(m-1)$, because $\varphi_m^{n+1}(m-1)=\varphi_m^n((m-1)
0)=\varphi_m^n(m-1)\varphi_m^n(0)$, and there is no letter in $A_m$ which is
smaller than~$0$.

For the other words in $SP_{L_m}$ when $n>2$, i.e., for $w_i$ with $0 \leq i <
{m-2}$, we follow Lemma~\ref{double_conv_seq_lem} by building a decreasing and
an increasing sequence of words converging to each of them.
Indeed, by construction of $w_i$ as $\varphi_m^\infty(i)$, we know that
$\varphi_m^{n}(i)$ is a prefix of $w_i$ for every $n \geq 0$.
Moreover, like above, $\varphi_m^{n}(i)$ is extendable to the right by $(i+1)$
in $w_i$ since $\varphi_m^{n+1}(i)= \varphi_m^{n}(i)\varphi_m^{n}(i+1)$, that
is, $\varphi_m^{n}(i)\:(i+1)$ is also a prefix of $w_i$.
Now, according to Lemma~\ref{bisp_words_lem}, $\varphi_m^{n}(i)$ is a bispecial
factor for every $n \geq 1$ such that $\#Biext=2m$,
and according to Lemma~\ref{sp_degree_lem} it must also be a
right special word with $\#Rext=m$.
Hence, $\varphi_m^{n}(i)\:(i+2)$ also belongs to $Fact_{L_m}$, and by definition
of $L_m$ there are words in $L_m$ whose prefix is $\varphi_m^{n}(i)\:(i+2)$.
Let us take one of them and denote it by $v_n$. 
In $\mathbb{Z}_m$, since $i < (m-2)$, we have that $(i+1) < (i+2)$, so that
$\varphi_m^{n}(i)\:(i+1) < \varphi_m^{n}(i)\:(i+2)$.
This means that $w_i < v_n$, and that $v_{n+1} < v_n$ for every $n \geq
1$.
Hence, $\{v_n\}_{n \in \mathbb{N}}$ is a decreasing sequence converging to $w_i$.

Similar arguments apply to build an increasing sequence converging to~$w_i$:
instead of $\varphi_m^{n}(i)\:(i+1)$, consider the factors $\varphi_m^{n}(i)\:i$,
for every $n \geq 1$, which also belong to $Fact_{L_m}$ and such that
$\varphi_m^{n}(i) \: i < \varphi_m^{n}(i) \: (i+1)$. Since $i \geq 0$, we have
that $i < (i+1)$ in $\mathbb{Z}_m$, and therefore, there exist words $v'_n$
in~$L_m$ for every $n \geq 1$ having $\varphi_m^{n}(i) \: i$ as prefixes, making
$\{v'_n\}_{n \in \mathbb{N}}$ as an increasing sequence converging to $w_i$.  \EPR

\vspace*{0.3cm}
\noindent
Now, putting together the obtained results so far, the finite non-injectivity of
$T_{L_m}$ can be described:

\begin{Theo}\label{non_inj_th}
  Let $L_m$ be a Thue-Morse subshift, $m\geq 2$, 
  and let $T_{L_m}$ be its semi-conjugate interval exchange transformation.
  Then $T_{L_m}$ is injective on $\II$ except on exactly a set of $m(m-1)$
  points which can be partitioned into $m-1$ subsets of $m$ points on which
  $T_{L_m}$ goes to the same image.
\end{Theo}
\PR According to Lemma~\ref{conseq_sp_lem}, $w_{m-2}$ and $w_{m-1}$ of
$SP_{L_m}$ are consecutive words,
so that by Lemma~\ref{consec_preserv_TM_lem}, $a_iw_{m-2}$ and
$a_iw_{m-1}$ in $LextSP_{L_m}$ must also be consecutive for every $a_i \in A_m$.
Hence, every pair $(a_iw_{m-2}, a_iw_{m-1}) \in \ConsecPairLm$, and we have that
only $a_iw_{m-1}$ belongs to $\Lmprim$ for every $i$.
But then, when $a_i \neq a_j$, $\phiLprim(\{a_iw_{m-1},a_jw_{m-1}\})$ is equal
to two distinct points.
We can then apply Lemma~\ref{noninj_lem} to conclude that non-injectivity of
$T_{L_m}$ occurs for these two points, i.e., they both belong to $\Ninjm$.
More generally, we obtain that $\phiLprim(Lext(\{w_{m-1}\}))$ is made of $m$
distinct points, all belonging in $\Ninjm$, on which $T_{L_m}$ goes to the same
image since $T_{L_m}(\phi_\mu(a_iw_{m-1})) = \phi_\mu(\sigma(a_iw_{m-1}))
= \phi_\mu(w_{m-1})$ for every $a_i \in A_m$.
This makes a first set of $m$ points out of the two words $w_{m-2}$ and
$w_{m-1}$ of $SP_{L_m}$ on which $T_{L_m}$ is non-injective.
Also, because of consecutivity, $w_{m-2}$ and the words in $Lext(w_{m-2})$ play
no role here as they do not belong to $\Lmprim$, being just the consecutive
smaller counterparts of respectively $w_{m-1}$ and the words in $Lext(w_{m-1})$.

For the other words $w_0,...,w_{m-3}$ of $SP_{L_m}$, according to
Lemma~\ref{conseq_sp_lem}, they are not part of any pair of consecutive words,
so that by Lemma~\ref{non-consec_preserv_TM_lem}, the same is true for all their
extensions $a_jw_i$ in $LextSP_{L_m}$, with $a_j \in A_m$.
Hence, all of them belong to $\Lmprim$.
We can then apply Lemma~\ref{noninj_lem} on $\phiLprim(Lext(w_{i}))$ to conclude
it is made of $m$ distinct points, all belonging in $\Ninjm$, on which $T_{L_m}$
goes to the same image, that is, $\phi_\mu(w_{i})$.
This makes the $m-2$ other sets of $m$ points on which $T_{L_m}$ is
non-injective.

Finally, all the above non-injectivity points are distinct since they all come
from distinct words in $\Lmprim$.  \EPR

\noindent
About these non-injectivity points of $T_{L_m}$, by
Lemma~\ref{ninj_discaccr_lem}, we also have:

\begin{Cor} 
Let $L_m$ be a Thue-Morse subshift, $m\geq 2$, 
and let $T_{L_m}$ be its semi-conjugate interval exchange transformation.
Then the set of the $m(m-1)$ non-injectivity points of $T_{L_m}$ described by
Theorem~\ref{non_inj_th} is exactly its set  $\Discaccr$ of discontinuities from
the right.
\end{Cor}
\section{Examples}
\label{ex_sec}

Given a subshift $L \in \ShiftsAMZ$, it is possible to obtain explicit
approximations of the $T_L$'s graph
by using $T_L$'s definition in Section~\ref{T_L_construct_sec}. In particular,
the following approximation technique works well when $L$ has been induced by 
primitive substitutions like the Thue-Morse ones:
Let $p_L(n)$ be the map giving the number of factors of length~$n$ occurring in
the words of~$L$,
i.e., $p_L(n)=\#Fact_L(n)$, called the {\em complexity} of~$L$~\cite{Fog02}.
For each $n>0$, let $T_L^{(n)}$ be the map over $\II$ whose graph is defined as
follows: its abscissa is divided into $p_{L}(n)$ right-open intervals of equal
length, and its ordinate is divided into $p_{L}(n-1)$ ones.
Following the order between factors, the intervals of the abscissa are put into
correspondence with the cylinders of $L$ determined by the factors in
$Fact_{L}(n)$, and the same is done for the intervals of the ordinate with the
cylinders of the factors in $Fact_{L}(n-1)$.
Then, $T_L^{(n)}$ is the piecewise affine map
which sends each interval corresponding to $Cyl_L(v)$, $v \in Fact_L(n)$ (i.e.,
$\phiLprim(Cyl_{L}(v))$, denoted $I_{v}$ from now), to the interval
corresponding to $Cyl_L(\sigma(v))$, i.e., $I_{\sigma(v)}$, using a slope
$\frac{p_{L}(n)}{p_{L}(n-1)}$.
Now, for subshifts $L$ like the Thue-Morse ones, with $p_L(n) = O(n)$%
~\cite{TS95},
uniquely ergodic~\cite{Que10},
the slopes $\frac{p_{L}(n)}{p_{L}(n-1)}$ can be proved to converge to~$1$ when
$n$ goes to infinity, and more generally that $T_L^{(n)}$ converges
to~$T_L$~\cite{LN17}.

Note that for every $n$, the map $T_L^{(n)}$ is surjective but not injective on
the intervals corresponding to factors of $Fact_{L}(n)$ in the form $au$ where
$u \in Fact_L(n-1)$ is left special.
However, such a factor $au$ growing with $n$ to the right is such that either it
eventually becomes $auu' \in Fact_L(n')$, $n'>n$, where $uu'$ is not left
special, so that $T_L^{(n')}$ becomes injective on $I_{auu'}$, or it converges
to a word of $LextSP_L$.
As expected then, $T_L^{(n)}$ becomes ``more and more injective'' as $n$
increases and as it converges to $T_L$.
And $T_L^{(n)}$ approximates well $T_L$ even for factors $au$ with a left
special~$u$ because sending $I_{au}$ to $I_{u}$ is close to what happens when
$au$ has been extended as $auu' \in Fact_L(n')$ with a non left special $uu'$
since $I_{auu'} \subseteq I_{au}$ and $I_{uu'} \subseteq I_{u}$.
Also, in order to visualize more faithfully the construction of~$T_L$ as given
in Section~\ref{T_L_construct_sec}, together with the accumulation points
corresponding to the words in $LextSP_L$, one can restrict $T_L^{(n)}$ to be
only defined on intervals corresponding to the cylinders of the partition
$PART_L$, based only on factors with non-left special suffixes, hence
corresponding to injectivity continuity intervals of $T_L$ and $T_L^{(n)}$.  We
denote by $T_{PART_L}^{(n)}$ this restriction of~$T_L^{(n)}$ based on the factors
determining $PART_L$ up to length~$n$
(see their iterative construction in p.~\pageref{part_L_algo}).
However, note that this faithfulness requires much larger $n$ than for
$T_L^{(n)}$ to be able to see images of $T_L$ for short injectivity intervals as
they generally correspond to very long factors determining $PART_L$.

Here is then what happens for the Thue-Morse $L_3$ case
based on the substitution: 
\[\begin{array}{llllll}
\varphi_3 (0) = 01, \;\;\; \varphi_3 (1) = 12, \;\;\;\varphi_3 (2) = 20.
\end{array}\]

\noindent
Its three infinite left special words $SP_{L_3}=Fix(\varphi_3)$ are: 
\[\begin{array}{llllll} 
w_0=0112122012202001122020012001011...,\\
w_1=1220200120010112200101120112122...,\\
w_2=2001011201121220011212201220200....
\end{array} \]

\noindent
According to Lemma~\ref{conseq_sp_lem}, $w_1$ and $w_2$ are consecutive words.
By Theorem~\ref{non_inj_th}, they induce three points of non-injectivity for
$T_{L_3}$ from $\phiLprim(Lext(\{w_1,w_2\}))$ such that (the extensions of $w_1$
taking no role as they do not belong to $\Lprim_3$, being consecutive
with their corresponding extensions of $w_2$):
\[\begin{array}{ll}
T_{L_3}(\phi_\mu(0w_2))=T_{L_3}(\phi_\mu(1w_2))=T_{L_3}(\phi_\mu(2w_2)).
\end{array}\]

\noindent
For $w_0$, the words in $Lext(w_0)$ are all in $\Lprim_3$ by
Lemma~\ref{conseq_sp_lem}, and Theorem~\ref{non_inj_th} says they also induce
three points of non-injectivity for $T_{L_3}$ from $\phiLprim(Lext(w_0))$:
\[\begin{array}{llllll}
T_{L_3}(\phi_\mu(0w_0))=T_{L_3}(\phi_\mu(1w_0))=T_{L_3}(\phi_\mu(2w_0)). 
\end{array}\]

\noindent
Then the function graph of $T_{L_3}$ looks like (as $T_{PART_{L_3}}^{(300)}$): 

\begin{center} 
  \includegraphics[width=11.5cm]{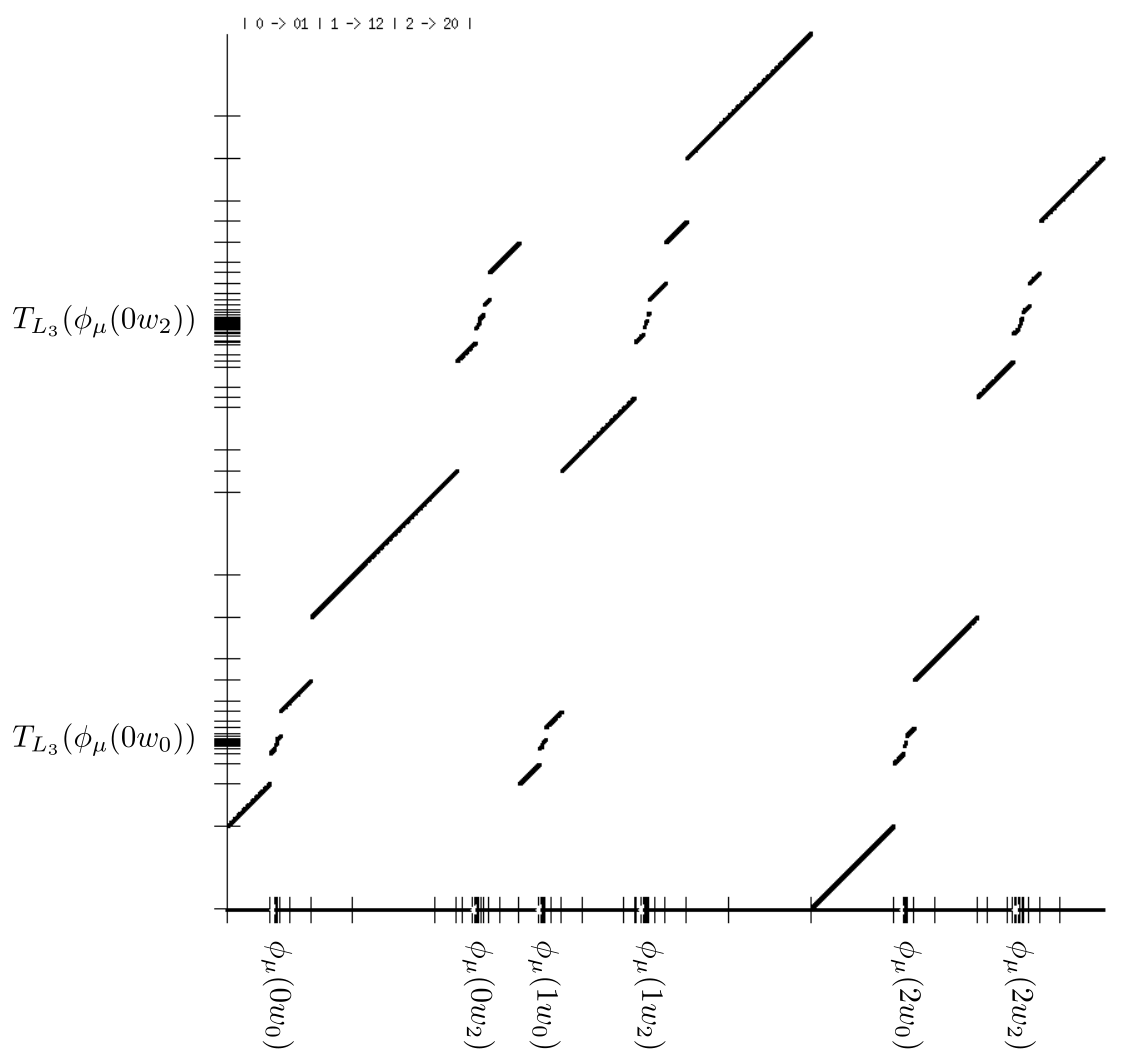}
  \label{thuegen3_300}
\end{center}

\noindent
In accordance with Theorem~\ref{non_inj_th}, its configuration of
non-injectivity is made of two sets of three points respectively going to
$T_{L_3}(\phi_\mu(0w_2)$ and $T_{L_3}(\phi_\mu(0w_0))$:

\begin{center} 
  \includegraphics[width=4cm]{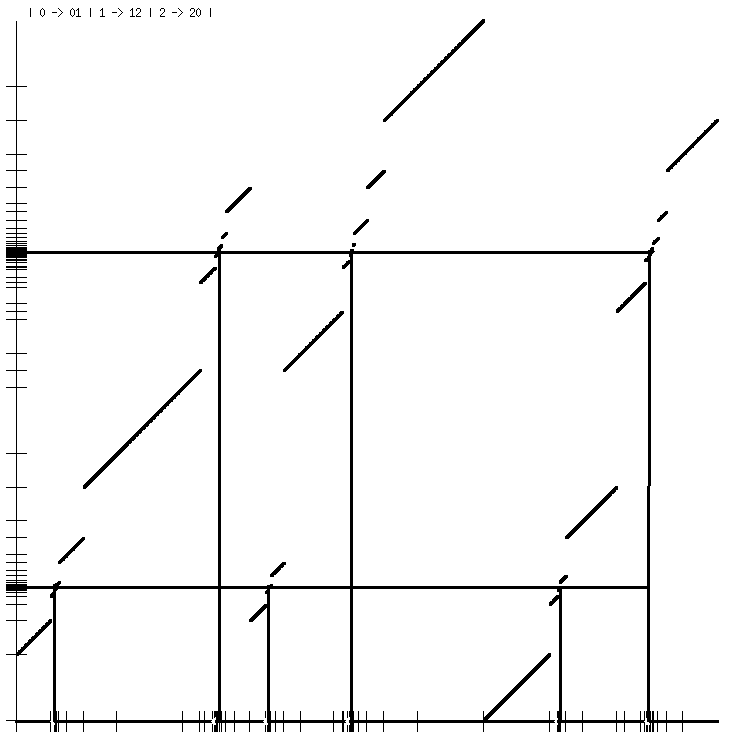}
\end{center}

\vspace*{0.3cm}
\noindent
Let us also see what happens for the $L_4$ case
based on the substitution: 
\[\begin{array}{llllll}
\varphi_4 (0) = 01, \;\;\; \varphi_4 (1) = 12, \;\;\; \varphi_4 (3) = 23,
\;\;\; \varphi_4 (2) = 30.
\end{array}\]

\noindent
The four infinite left special words $SP_{L_4}=Fix(\varphi_4)$ of $L_4$ are:
\[\begin{array}{llllll}
w_0=0112122312232330122323302330300...\\
w_1=1223233023303001233030013001011...\\
w_2=2330300130010112300101120112122...\\
w_3=3001011201121223011212231223233...
\end{array}\]

\noindent
According to Lemma~\ref{conseq_sp_lem}, $w_2$ and $w_3$ are consecutive words,
and by Theorem~\ref{non_inj_th}, they induce four points of non-injectivity for
$T_{L_4}$ such that:
\[\begin{array}{llllll}
T_{L_4}(\phi_\mu(0w_3))= T_{L_4}(\phi_\mu(1w_3))= T_{L_4}(\phi_\mu(2w_3))=
T_{L_4}(\phi_\mu(3w_3)).
\end{array}\]

\noindent
For $w_0$ and $w_1$, each set of extensions $Lext(w_0)$ and $Lext(w_1)$ induce
four points of non-injectivity for $T_{L_4}$:
\[\begin{array}{ll}
T_{L_4}(\phi_\mu(0w_0))=T_{L_4}(\phi_\mu(1w_0))=T_{L_4}(\phi_\mu(2w_0))=T_{L_4}(\phi_\mu(3w_0)),\\
T_{L_4}(\phi_\mu(0w_1))=T_{L_4}(\phi_\mu(1w_1))=T_{L_4}(\phi_\mu(2w_1))=T_{L_4}(\phi_\mu(3w_1)). 
\end{array}\]

\noindent
Then the function graph of $T_{L_4}$ looks like (as $T_{PART_{L_4}}^{(300)}$):

\begin{center} 
  \includegraphics[width=11.5cm]{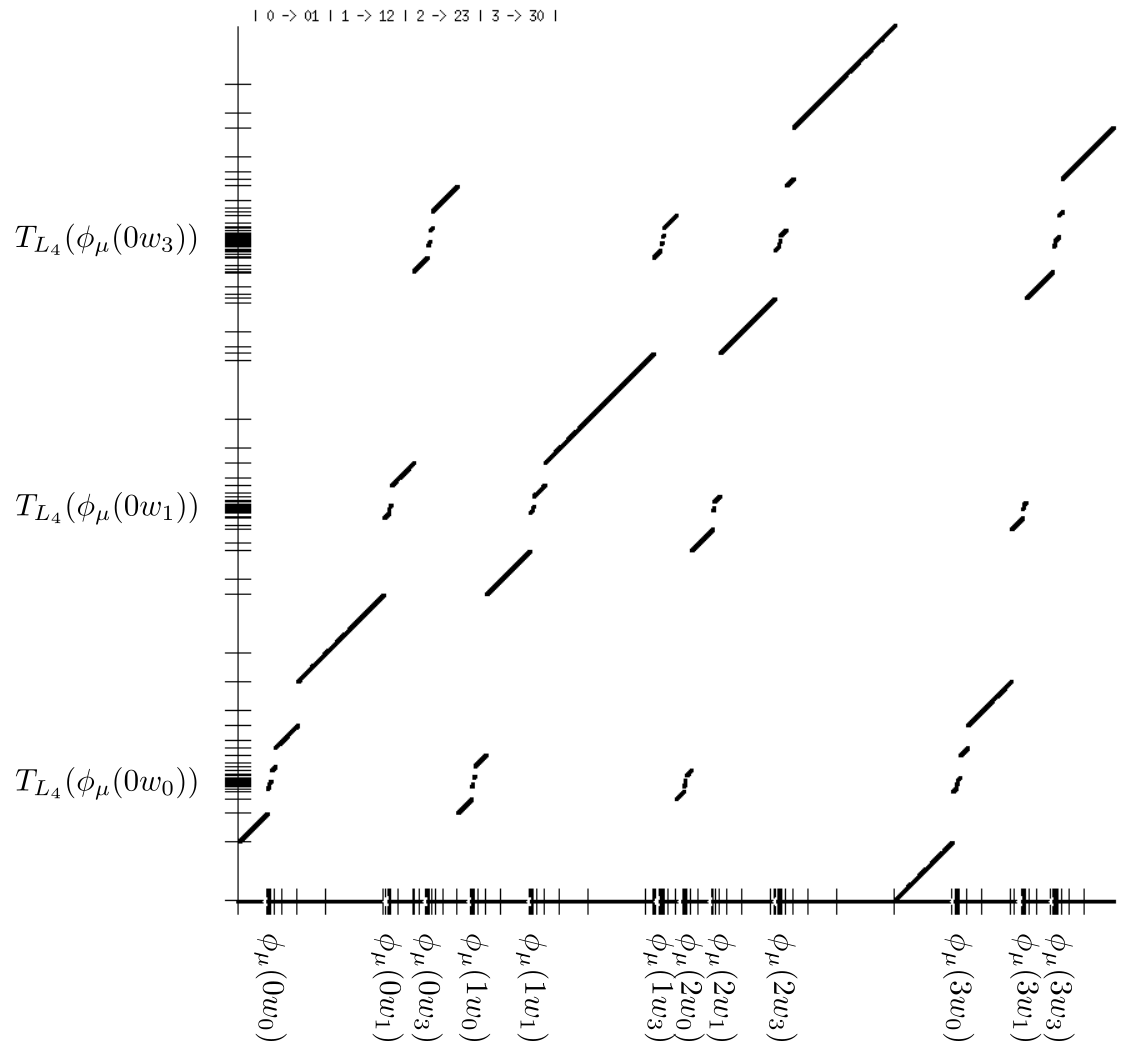}
  \label{thuegen4_300}
\end{center}

\noindent
Accordingly, its configuration of non-injectivity is
made of three sets of four points respectively going to
$T_{L_4}(\phi_\mu(0w_3))$, $T_{L_4}(\phi_\mu(0w_0))$ and
$T_{L_4}(\phi_\mu(0w_1))$:
\begin{center} 
  \includegraphics[width=4cm]{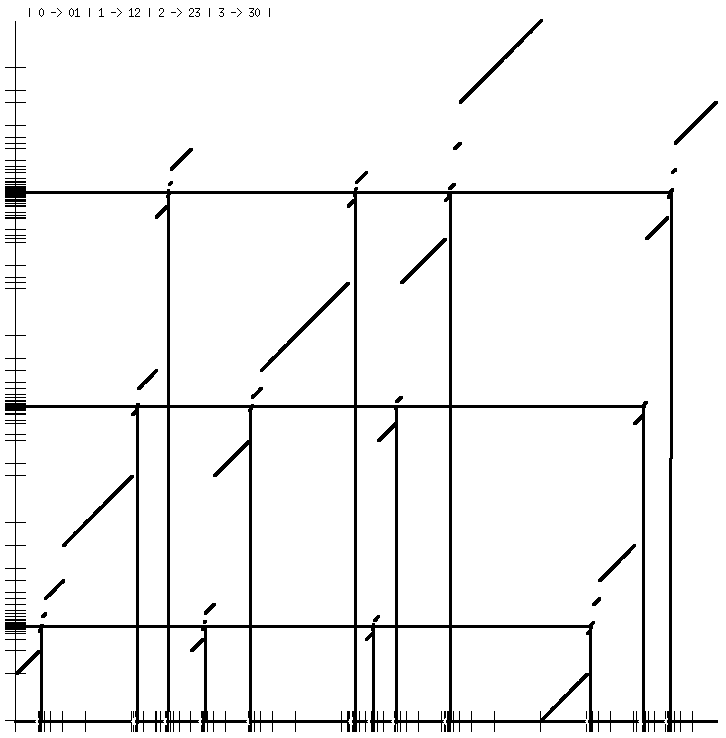}
\end{center}

\vspace*{0.3cm}
\noindent
Going further with larger $m$'s, here is how looks the function graph of
$T_{L_8}$ (as $T_{PART_{L_8}}^{(200)}$)
whose configuration of non-injectivity is made of seven sets of eight points
going to the same image:

\begin{center} 
  \includegraphics[width=10cm]{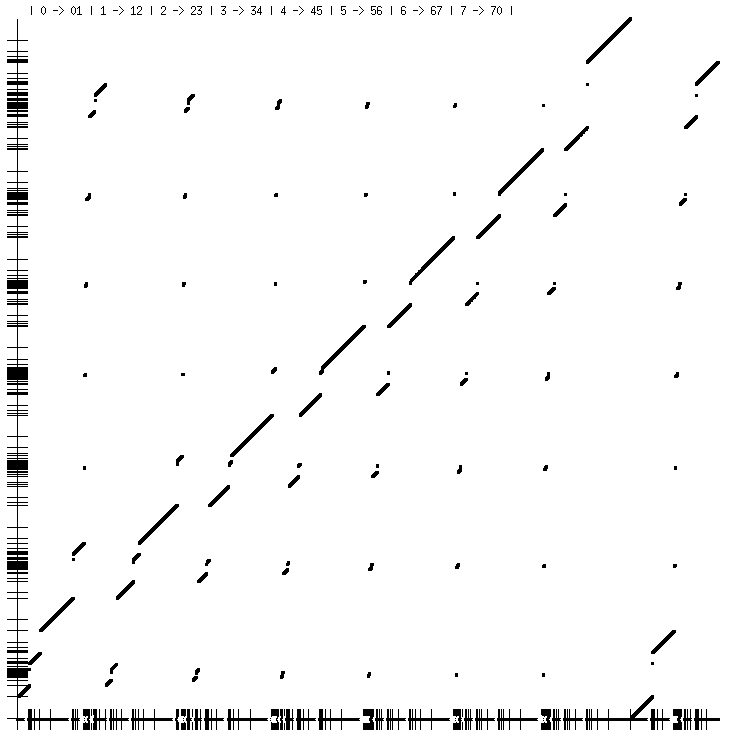}
\end{center}

\vspace*{0.3cm}

\noindent
Here are also how respectively look the graphs of the interval exchange
transformations $T_{L_{15}}$ and $T_{L_{30}}$ (as $T_{L_{15}}^{(50)}$ and
$T_{L_{30}}^{(50)}$ --~using $T_{PART_{L_{.}}}^{(n)}$ here would require much
larger $n$ without much gain because of the high density of accumulation points
and the shortness of most of the injectivity intervals):

\vspace*{0.2cm}
\begin{center} 
  \includegraphics[width=6.1cm]{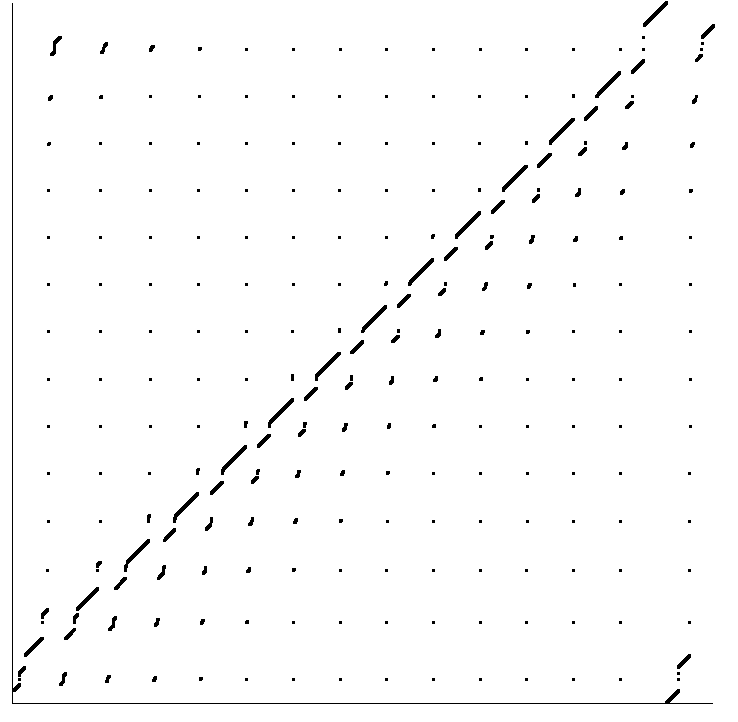}
  \hspace*{1cm}
  \includegraphics[width=6.1cm]{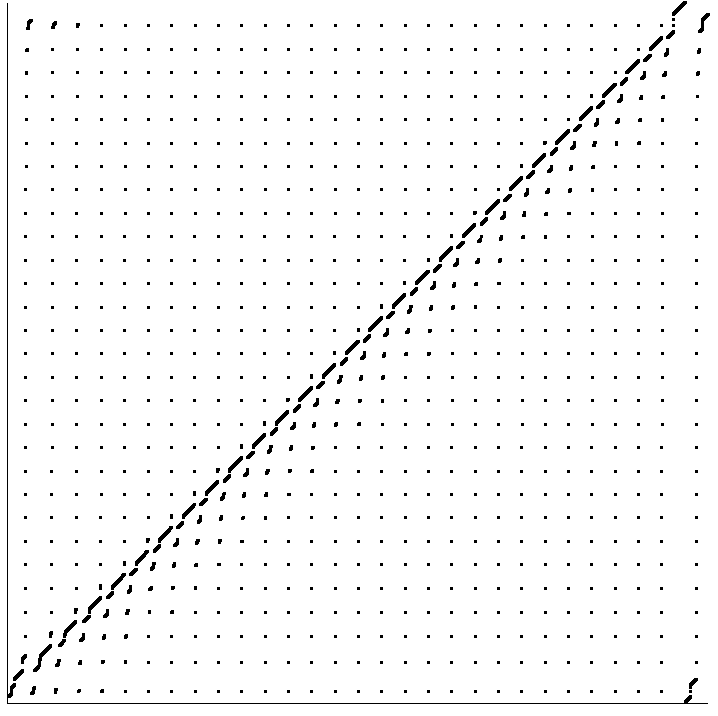}
\end{center}

{\footnotesize
\providecommand{\bysame}{\leavevmode\hbox to3em{\hrulefill}\thinspace}
\providecommand{\MR}{\relax\ifhmode\unskip\space\fi MR }
\providecommand{\MRhref}[2]{%
  \href{http://www.ams.org/mathscinet-getitem?mr=#1}{#2}
}
\providecommand{\href}[2]{#2}

}

\end{document}